\documentclass[12pt,psamsfonts]{amsart}

\newtheorem{anyprop}{Anyprop}[section]

\newtheorem{theorem}[anyprop]{Theorem}
\newtheorem{lemma}[anyprop]{Lemma}

\newtheorem{corollary}[anyprop]{Corollary}

\theoremstyle{definition}

\newtheorem{definition}[anyprop]{Definition}

\newtheorem{remark}[anyprop]{Remark}

\newcommand{\NN}{\mathbb{N}}

\newcommand{\FF}{\mathbb{F}}

\newcommand{\PP}{\mathbb{P}}

\newcommand  {\shF}     {\mathcal{F}}
\newcommand  {\shG}     {\mathcal{G}}

\newcommand  {\shK}     {\mathcal{K}}

\newcommand  {\shL}     {\mathcal{L}}

\newcommand  {\shS}     {\mathcal{S}}
\newcommand  {\shT}     {\mathcal{T}}

\newcommand  {\Char}    {\operatorname{char}}

\newcommand  {\dual}    {\vee}

\newcommand  {\lra}     {\longrightarrow}

\newcommand  {\modu}     {\operatorname{mod}}

\renewcommand{\O}       {\mathcal{O}}

\newcommand  {\Proj}    {\operatorname{Proj}}

\newcommand  {\ra}      {\rightarrow}

\newcommand  {\rk}    {\operatorname{rk}}

\newcommand  {\Syz}     {\operatorname{Syz}}

\newcommand{\numiii}{\renewcommand{\labelenumi}{(\roman{enumi})}}

\newcommand{\comdots}{ , \ldots , }
\newcommand{\komdots}{ , \ldots , }
\newcommand{\plusdots}{ + \ldots + }

\newcommand{\oplusdots}{ \oplus \ldots \oplus }

\newcommand{\subsetdots}{ \subset \ldots \subset }

\theoremstyle{remark}

\numberwithin{equation}{section}

\usepackage{amscd}
\usepackage{amssymb}

\def\mydate{\number\day\space\ifcase\month \or January\or February\or March\or April\or May\or
June\or July\or August\or September\or October\or November\or
December\fi \space\number\year}

\newcommand{\sheS} {\shT}

\newcommand{\boundo}{k}
\newcommand{\boundb}{b}
\newcommand{\bounda}{b}
\newcommand{\boundc}{k}
\newcommand{\boundr}{k}
\newcommand{\boundhn}{h}
\newcommand{\boundf}{k}

\newcommand{\boundm}{m}
\newcommand{\boundtight}{b}
\newcommand{\boundn}{n}
\newcommand{\bonum}{u}
\newcommand{\cu}{C}
\newcommand{\frob}{\Phi}
\newcommand{\fifi}{\FF}
\newcommand{\test}{\tau}
\newcommand{\minpricomp}{\circ}
\newcommand{\nobound}{b}

\begin{document}

\title[Bounds for test exponents]
{Bounds for test exponents}

\author[Holger Brenner]{Holger Brenner}
\address{Department of Pure Mathematics, University of Sheffield,
  Hicks Building, Hounsfield Road, Sheffield S3 7RH, United Kingdom}
\email{H.Brenner@sheffield.ac.uk}


\date{\mydate}



\begin{abstract}
Suppose that $R$ is a two-dimensional normal standard-graded domain
over a finite field. We prove that there exists a uniform Frobenius
test exponent $\boundb$ for the class of homogeneous ideals in $R$
generated by at most $n$ elements. This means that for every ideal
$I$ in this class we have that $f^{p^\boundb} \in I^{[p^\boundb]}$
if and only if $f \in I^F$. This gives in particular a finite test
for the Frobenius closure. On the other hand we show that there is
no uniform bound for Frobenius test exponent for all homogeneous
ideals independent of the number of generators. Under similar
assumptions we prove also the existence of a bound for tight closure
test ideal exponents for ideals generated by at most $n$ elements.
\end{abstract}

\maketitle



\section*{Introduction}

In this paper we deal with test exponents for Frobenius closure and
for tight closure. In this introduction we will concentrate on the
Frobenius closure and we will come back to tight closure in section
\ref{tightexponent}. Suppose that $R$ denotes a commutative
noetherian ring containing a field of positive characteristic $p$,
and let $I \subseteq R$ denote an ideal. The Frobenius closure,
written $I^F$, is the ideal
$$I^F=\{ f \in R: \, \exists \,  e \in \NN \,
\mbox{ such that } \, f^{p^{e}} \in I^{[p^{e}]} \}\, .$$ Here the
Frobenius power $I^{[q]}=(f^q: f \in I)$ is the extended ideal of
$I$ under the $e$-th iteration of the Frobenius homomorphism $R \ra
R$, $f \mapsto f^q$, where $q=p^{e}$. Since $I^F$ is an ideal in a
noetherian ring, there exists a number $b$ such that $f^{p^b} \in
I^{[p^b]}$ for every $f \in I^F$.

A problem of Katzman and Sharp
\cite[Introduction]{katzmansharpfrobenius} asks in its strongest
form: does there exist a number $b$, depending only on the ring $R$,
such that, for every ideal and for every $f \in I^F$, we have
$f^{p^b} \in I^{[p^b]}$. A positive answer to this question,
together with the actual knowledge of a bound for $b$, would give an
algorithm to compute the Frobenius closure $I^F$. We call such a
number $b$ a Frobenius test exponent for the ring $R$.

A weaker question is whether for a given ideal $I$ there exists a
Frobenius test exponent for all Frobenius powers $I^{[q]}$. The
existence of such a weak Frobenius test exponent $b$ means that $f
\in (I^{[q]})^F$ holds if and only if $f^{p^b} \in
(I^{[q]})^{[p^b]}$. Katzman and Sharp show that for an ideal
generated by a regular sequence there exists such a weak Frobenius
test exponent.

An intermediate question is the following: does there exist a
Frobenius test exponent for ideals generated by at most a fixed
number of elements. Since $(f_1 \komdots f_n)^{[q]}=(f_1^q \komdots
f_n^q)$, the Frobenius powers of an ideal generated by $n$ elements
do not need more generators than $n$; therefore a positive answer to
this question implies a positive answer to the weak question of
Katzman and Sharp. Moreover, a positive answer to this question has
the same computational impact on the decision whether $f \in I^F$
holds or not.

In this paper we consider the case of homogeneous $R_+$-primary
ideals $I$ in a geometrically normal two-dimensional standard-graded
domain $R$ over a field $\fifi$ of positive characteristic. We
obtain the following results. If $R$ is the homogeneous coordinate
ring of an elliptic curve, then the answer to the third question is
affirmative: there is a Frobenius test exponent $b$ for the class of
homogeneous $R_+$-primary ideals generated by at most $n$ elements,
and in fact one can take $b=n-1$ (Theorem \ref{testelliptic}). This
is an application of the Atiyah classification of vector bundles
over an elliptic curve.

In Theorem \ref{nofrobeniusbound} we give a negative answer to the
strong question of Katzman and Sharp. We show that for a homogeneous
coordinate ring $R$ over an elliptic curve with Hasse-invariant $0$
there does not exist a uniform Frobenius test exponent valid for all
homogeneous ideals in $R$. This relies on results of Oda about
properties of the Frobenius pull-backs of the bundles $F_r$ on an
elliptic curve, where $F_r$ denotes the unique indecomposable bundle
of rank $r$ and degree $0$ with a non-trivial section, and on
realizing $F_r$ as a syzygy bundle for some ideal generators.

For homogeneous normal coordinate rings of curves of higher genus we
obtain in Theorem \ref{testhighergenus} a positive answer to the
third question under the condition that the base field is finite.
This is deduced from the geometric interpretation of tight closure
in terms of vector bundles on the corresponding smooth projective
curve and the boundedness of semistable bundles on the curve of
given rank and degree. The boundedness implies that there is a
finite test for strong semistability if everything is defined over a
finite field (Lemma \ref{strongtest}). From this we deduce using a
recent result of V. Trivedi that there is also a finite algorithm to
find the strong Harder-Narasimhan filtration of a bundle of given
rank (Corollary \ref{hnbound}).

These two results are important ingredients in an argument to show
that a cohomology class $c \in H^1(\cu, \shS)$ of a locally free
sheaf of rank $r$ on a smooth projective curve defined over a finite
field which is annihilated by some Frobenius power is in fact
annihilated by a certain fixed power which depends only on the rank
and on the curve, but not on the sheaf itself (Theorem
\ref{testexponentcurve}). This yields eventually the Frobenius test
exponent for ideals generated by at most $n$ elements (Theorem
\ref{testhighergenus}). However, this Frobenius test exponent is --
even for ideals generated by three elements in the coordinate ring
over a curve of genus $2$ -- hardly a basis for implementing an
algorithm to compute the Frobenius closure.

In the same spirit we prove in Theorem \ref{tightexponenttheorem}
the existence of a uniform test ideal exponent for ideals generated
by $\leq n$ elements, under the assumption that $R$ is a
two-dimensional standard-graded geometrically normal Gorenstein
domain over a finite field.

I would like to thank M. Hochster, M. Katzman, R. Sharp, K. Smith,
the referee and in particular R. Lazarsfeld for their interest and
for useful comments.

\section{Frobenius test exponents over elliptic curves}

The vector bundles (locally free sheaves) on an elliptic curve $\cu$
over an algebraically closed field are well understood due to the
classification of Atiyah; see \cite{atiyahelliptic}. We recall
briefly some consequences of this classification. If $\shS$ is an
indecomposable bundle, meaning that there is no non-trivial
decomposition $\shS= \shF \oplus \shG$ with vector bundles $\shF$
and $\shG$ of smaller ranks, then the following hold: if $\deg
(\shS)
>0$, then $H^1(\cu,\shS)=0$. If $\deg (\shS) =0$, then either
$H^0(\cu,\shS)=H^1(\cu,\shS)=0$ or $\shS=F_r$, the unique
indecomposable sheaf of rank $r$ and degree $0$ with a non-trivial
global section. These sheaves have $\dim H^0 (\cu, F_r)= \dim
H^1(\cu,F_r)=1$, they are self-dual and they are related by the
non-trivial extensions $0 \ra F_{r-1} \ra F_r \ra \O_{\cu} \ra 0$.

\begin{theorem}
\label{testelliptic} Let $R$ denote a normal standard-graded
two-dimensional domain over a field $\fifi$ of positive
characteristic $p$ and suppose that $\cu= \Proj R$ is an elliptic
curve. Let $I=(f_1 \komdots f_n)$ denote an $R_+$-primary ideal
generated by $n$ homogeneous elements. Suppose that the element $f$
belongs to the Frobenius closure of $I$. Then already $f^q \in
I^{[q]}$ holds for $q= p^{n-1}$. Hence $n-1$ is a uniform Frobenius
test exponent for all homogeneous $R_+$-primary ideals generated by
$\leq n$ elements.
\end{theorem}
\proof We may assume that $\fifi$ is algebraically closed. Since $I$
is homogeneous, all Frobenius powers $I^{[q]}$ and the Frobenius
closure $I^F$ are also homogeneous. Therefore we may assume that $f$
is homogeneous. A homogeneous element $f$ of degree $m$ yields via
the short exact sequence
$$0 \lra \Syz(f_1 \komdots f_n)(m) \lra \bigoplus_{i=1}^n
\O_{\cu}(m-d_i) \stackrel{f_1 \komdots f_n}{\lra} \O_{\cu}(m) \lra 0
$$ of locally free sheaves on $\cu$ the class $\delta(f)
\in H^1(\cu, \Syz(f_1 \komdots f_n)(m))$. The pull-back of this
sequence under the $e$-th absolute Frobenius morphism $\frob^{e}: C
\ra C$ is ($q=p^{e}$)
$$0 \lra \frob^{e*}(\Syz(f_1 \komdots f_n)(m)) \lra \bigoplus_{i=1}^n \O_{\cu}(qm-qd_i)
\stackrel{f_1^q \komdots f_n^q}{\lra} \O_{\cu}(qm) \lra 0 \, ,$$
therefore $\Syz(f_1^q \komdots f_n^q)(qm) =\frob^{e*}(\Syz(f_1
\komdots f_n)(m))$. From this we see that $f^q \in I^{[q]}=(f_1^q
\komdots f_n^q)$ if and only if the $e$-th Frobenius pull-back of
the class $\delta(f)$ is zero, that is $\frob^{e*}(\delta(f)) =0$ in
$H^1(\cu, \Syz(f_1^q \komdots f_n^q)(qm))$. The element $f$ belongs
to the Frobenius closure of $I=(f_1 \komdots f_n)$ if and only if
some Frobenius power of this cohomology class vanishes.

Let $\Syz(f_1 \komdots f_n)(m)= \shS_1 \oplusdots \shS_k$ denote the
decomposition of the syzygy bundle into indecomposable bundles, and
denote the components of the cohomology class $c=\delta(f)$ by
$c_j\in H^1(\cu,\shS_j)$. For the $\shS_j$ of positive degree we
have $H^1(\cu, \shS_j)=0$. If $\shS_j$ has negative degree and $c_j
\neq 0$, then $f$ does not belong to the tight closure of the ideal
by \cite[Corollary 4.1]{brennertightelliptic} and so it does not
belong to the Frobenius closure. So these components are also zero.
For an indecomposable sheaf $\shS$ of degree $0$ we have two
possibilities depending on whether $\shS$ has no non-trivial global
sections or it has. In the first case we have again
$H^1(\cu,\shS)=0$. In the second case we have $\shS =F_r$, the
unique indecomposable sheaf of rank $r$ and degree $0$ with a
non-trivial global section. Furthermore $r \leq n-1 = \rk (\Syz(f_1
\komdots f_n)(m))$.

The indecomposable bundles $F_r$ are related by short exact
sequences $$0 \lra F_{r-1} \lra F_r \lra \O_{\cu} \lra 0 \, ,$$
where $1 \in \fifi = \Gamma(\cu, \O_{\cu})$ maps to a non-trivial
element in the one-dimensional space $H^1(\cu, F_{r-1})$. Suppose
first that the Hasse-invariant of the elliptic curve is $0$, so that
the Frobenius morphism $\frob^* : H^1(\cu, \O_{\cu}) \ra H^1(\cu,
\O_{\cu})$ is the zero map. In this case we show by induction on $r$
that the $(r-1)$-th Frobenius pull-back trivializes $F_r$ (meaning
that $\frob^{(r-1)*} (F_r) \cong \O_{\cu}^r$) and that the $r$-th
Frobenius pull-back annihilates $H^1(\cu, F_r)$. For $F_1 =\O_{\cu}$
this follows from the assumption that the Hasse-invariant is $0$.
Now suppose that $r  \geq 2$. The $(r-2)$-th Frobenius pull-back of
$0 \ra F_{r-1} \ra F_r \ra \O_{\cu} \ra 0$ yields an exact sequence
$$0 \lra \O_{\cu}^{r-1} \lra \frob^{(r-2)*} (F_r) \lra \O_{\cu} \lra 0 \, . $$
This extension is given by a cohomology class in $H^1(\cu,
\O_{\cu}^{r-1})$, and so application of the Frobenius once more
shows that the extension is trivial; hence $\frob^{(r-1)*} (F_r)
\cong \O_{\cu}^r$. Another application of the Frobenius annihilates
then the cohomology coming from $F_r$.

If the Hasse-invariant of $\cu$ is $1$, then the Frobenius acts
bijectively on $H^1(\cu, \O_{\cu})$. The short exact sequences
relating $F_{r-1}$ and $F_r$ induce isomorphisms $H^1(\cu, F_r)
\cong H^1(\cu, \O_{\cu})$. Hence by induction on $r$ we get that
$\frob^*(F_r) \cong F_r$ and that $\frob^*: H^1(\cu, F_r) \ra
H^1(\cu, \frob^*(F_r))$ is a bijection. Hence the Frobenius does not
annihilate anything. \qed

\begin{remark}
Theorem \ref{testelliptic} implies that the Frobenius closure of a
homogeneous $R_+$-primary ideal in a normal homogeneous coordinate
ring $R$ of an elliptic curve can be computed by an easy algorithm.
The same argument shows (for Hasse-invariant $0$ or $1$) that an
element $f$ belongs to the plus closure $I^+$ (or to the tight
closure) of $I$ if and only if it belongs to the extended ideal
under the $(n-1)$-th iteration of the ring homomorphism $\varphi: R
\ra R$ which describes the $p$-multiplication of the elliptic curve.
This shows that the tight closure of an $R_+$-primary homogeneous
ideal in $R$ is also algorithmically computable. Another easy
algorithm for the computation of tight closure in
$R=\fifi[x,y,z]/(G)$, where $G$ is a cubic polynomial, in terms of
test ideal exponent is given in Corollary \ref{testcubic} below.
\end{remark}

\section{Unboundedness of Frobenius test exponents}

In this section we show that in a homogeneous coordinate ring $R$ of
an elliptic curve with Hasse-invariant $0$ there does not exist a
bound for Frobenius test exponents which holds uniformly for all
homogeneous ideals.

We recall the following result of Oda about the Frobenius pull-backs
of the bundles $F_r$ on an elliptic curve and how the Frobenius acts
on the cohomology \cite[Proposition 2.10, Theorem
2.17]{odaelliptic}.

\begin{theorem}
\label{oda} Let $\cu$ denote an elliptic curve over a field $\fifi$
of positive characteristic $p$ and let $F_r$ denote the unique
indecomposable bundle on $\cu$ of rank $r$ and degree $0$ with
$H^0(\cu,F_r) =\fifi \neq 0$. Let $\frob: \cu \ra \cu $ denote the
Frobenius. Then the following hold.

\numiii
\begin{enumerate}

\item
If the Hasse-invariant of $\cu$ is non-zero, then $\frob^*(F_r)
\cong F_r$.

\item
If the Hasse-invariant of $\cu$ is zero, then
$$\frob^*(F_r) \cong\O_\cu^r \,\,\, \mbox { for }\,\,\, r \leq p \,\,\,\mbox{ and }\,\,\,
\frob^*(F_r) \cong \bigoplus_{i=1}^p F_{ \lfloor \frac{r-i}{p}
\rfloor +1} \,\,\, \mbox { for }\,\,\, r > p \, .$$

\item
The map $\frob^*: H^1(\cu,F_r) \ra H^1(\cu, \frob^*(F_r))$ is
injective if and only if the Hasse-invariant is non-zero or $r \geq
p$ {\rm(}and for Hasse-invariant zero and $r <p$ it is the zero
map{\rm)}.
\end{enumerate}
\end{theorem}

From this result of Oda we deduce the corollary that there does not
exist a uniform bound for Frobenius test exponents for all locally
free sheaves on $C$ independent of the rank.

\begin{corollary}
\label{odacor}
Let $\cu$ denote an elliptic curve over a field
$\fifi$ of positive characteristic $p$ and suppose that the
Hasse-invariant is $0$. Then for every number $\nobound \in \NN$
there exists a locally free sheaf $\shS$ on $\cu$ and a cohomology
class $c \in H^1(\cu,\shS)$ such that $c$ is annihilated by some
Frobenius power, but such that $\frob^{\nobound*}(c) \neq 0$.
\end{corollary}
\proof We show that for $0 \neq c \in H^1(\cu,F_r)$, where $r \geq
p^{\nobound+1}$, we have that $\frob^{\nobound*}(c) \in H^1(\cu,
\frob^{\nobound*} (F_r))$ is not zero, though it is annihilated by
some Frobenius power, since the curve is assumed to have
Hasse-invariant $0$ (see the proof of Theorem \ref{testelliptic}).
We use induction on $\nobound$, the case where $\nobound =0$ being
clear. So suppose that $\nobound \geq 1$. By Theorem
\ref{oda}(ii),(iii) we have the decomposition
$$\frob^*(F_r) \cong\bigoplus_{i=1}^p F_{ \lfloor \frac{r-i}{p}\rfloor +1}$$
and we know that the mapping on the cohomology is injective.  Hence
at least one component of $\frob^*(c)$ is non-zero, say $c' =
(\frob^*(c))_i \in H^1(\cu, F_{ \lfloor \frac{r-i}{p}\rfloor +1})$.
The rank of $F_{ \lfloor \frac{r-i}{p}\rfloor +1}$ is at least
$$ \big\lfloor \frac{r-i}{p} \big\rfloor +1 \geq \big \lfloor \frac{r-p}{p} \big\rfloor +1
=\big \lfloor \frac{r}{p}\big \rfloor \geq \big \lfloor
\frac{p^{\nobound+1}}{p} \big \rfloor = p^\nobound \, .$$ By the
induction hypothesis we know that $\frob^{ (\nobound-1)*}(c')\neq
0$; hence $ \frob^{ \nobound*} (c) = \frob^{(\nobound
-1)*}(\frob^*(c))\neq 0$. \qed

\medskip
We are going to translate these results into results about the
Frobenius closure of ideals. To do so we have to realize the bundles
$F_r$ on an elliptic curve as syzygy bundle for suitable ideal
generators. The following lemma is known.

\begin{lemma}
\label{bundleassyz} Let $C$ denote a smooth projective curve over an
algebraically closed field, and let $\shS$ denote a locally free
sheaf of rank $r$, which is globally generated by $r+k$ global
sections, $k \geq 1$. Then it is also globally generated by $r+1$
global sections.
\end{lemma}
\proof The assumption means that we have a surjection
$$ \O_C^{r+k} \lra \shS \lra 0   \, .$$
Set $V= \Gamma(C, \O_C^{r+k})$, and consider the family of subspaces
$U \subset V$ of codimension one, which form an $r+k$-dimensional
Grassmann variety $G(r+k-1,V)$. Fix $x \in C$ and consider the
surjection
$$\varphi_x : V \lra W := \shS \otimes \kappa(x) \, .$$
Let $K=K_x $ denote its kernel. The intersection of $U$ with $K$ has
either dimension $k-1$ or it is $K$. In the first case $U$ maps onto
$W$, in the second case not. The subspaces $U$ of dimension $r+k-1$
which contain $K$ are given by the $(r-1)$-dimensional subspaces of
$V/K$, so they form an $r$-dimensional subvariety $G_x$ of
$G(r+k-1,V)$. Hence $\bigcup_{x \in C} G_x$ has dimension at most
$r+1$, and so for $k \geq 2$ the generic hyperplane $U \subset V$
also generates $\shS$ globally. So inductively we can reduce $k$
until we have only $r+1$ generators. \qed

\begin{theorem}
\label{nofrobeniusbound} Let $R$ denote a normal homogeneous
coordinate ring over an elliptic curve of positive characteristic
$p$ and of Hasse-invariant $0$. Then there does not exist a
Frobenius test exponent valid for all homogeneous $R_+$-primary
ideals $I \subseteq R$.
\end{theorem}
\proof Let $C=\Proj R$ denote the corresponding elliptic curve and
let $\O_C(1)$ denote the ample invertible sheaf on $C$. For some
$\ell
>0$ the twisted bundle $F_r(\ell)$ on $C$ is globally generated. Due to Lemma
\ref{bundleassyz} there exists a short exact sequence
$$0 \lra \O_C(-d) \lra \O_C^{r+1} \lra F_r(\ell) \lra 0 \,  \,\,  (d  >0). $$
Dualizing and tensoring with $\O(\ell)$ yields (since $F_r^\dual
=F_r$)
$$0 \lra F_r \lra \O_C^{r+1} (\ell) \lra \O_C(d+ \ell) \lra 0 \,   ,$$
where the last mapping is given by some homogeneous elements $f_1
\comdots f_{r+1}$ of degree d, hence $F_r \cong \Syz(f_1 \komdots
f_{r+1})(d+\ell)$. Since $\ell >0$, we have $H^1(C, \O(\ell))=0$ and
therefore the non-trivial class $c \in H^1(C, F_r)$ comes from an
element $f \in \Gamma(C, \O_C(d+\ell))=R_{d+ \ell}$. By Corollary
\ref{odacor} and its proof there exists for given $b$ a number $r$
such that $0 \neq c \in H^1(C, F_r)$ is not annihilated by the
$b$-th Frobenius power, but it is annihilated by some Frobenius
power. This means that $f$ belongs to the Frobenius closure of the
ideal $I=(f_1 \komdots f_{r+1})$, but $f^{p^b} \not\in I^{[p^b]}$.
\qed

\section{Semistable sheaves on curves over a finite field}

We look now at normal homogeneous coordinate rings over a smooth
projective curve $\cu$ of higher genus. The theory of vector bundles
on $\cu$ is still the main ingredient in the following results, but
since this theory is more complicated than in the elliptic case we
obtain our results only under the condition that everything is
defined over a finite field.

Recall some facts (see  \cite{huybrechtslehn} or
\cite{seshadrifibre} for details) about locally free sheaves on a
smooth projective curve $\cu$ over a field $\fifi$. The degree of a
locally free sheaf $\shS$ on $\cu$ of rank $r$ is defined by $\deg
(\shS) = \deg \bigwedge^r (\shS)$; the degree is additive on short
exact sequences. The slope of $\shS$, written $\mu(\shS)$, is
defined to be $\deg(\shS)/ r$. The slope has the property that
$\mu(\shS \otimes \shT)= \mu(\shS) + \mu(\shT)$.

A locally free sheaf $\shS$ is called semistable (respectively
stable) if $\mu(\shT) \leq \mu(\shS)$ (respectively
$\mu(\shT)<\mu(\shS)$) for every locally free subsheaf $\shT \subset
\shS$. Tensoring with an invertible sheaf does not affect this
property. A non-zero morphism $\shS \ra \shT$ between two semistable
sheaves enforces $\mu(\shT) \geq \mu(\shS)$. In particular, a
semistable sheaf of negative degree does not have any non-trivial
global section.

For fixed rank $r$ and degree $d$ the set of semistable sheaves form
a bounded family: see \cite[1.7]{huybrechtslehn} or \cite[III.
A]{seshadrifibre}. This is a basic result in the construction of the
moduli space of (semi)stable sheaves on a curve and on higher
dimensional varieties. This boundedness implies in particular that
there are only finitely many semistable sheaves of rank $r$ and of
degree $d$ defined over a finite field.


In positive characteristic $p$, a locally free sheaf $\shS$ is
called strongly semistable if every Frobenius pull-back
$\frob^{e*}(\shS)$ under the $e$-th absolute Frobenius $\frob: \cu
\ra \cu$ is again semistable. The following Lemma shows that if
everything is defined over a finite field, then we only have to look
at certain Frobenius powers to test for strong semistability.

\begin{lemma}
\label{strongtest} Let $\cu$ denote a smooth projective curve
defined over a finite field $\fifi$, and let $r \in \NN$. Then there
exists a number $\boundm$ such that for every locally free sheaf
$\shS$ of rank $\, r$ on $\cu$ the following holds: if $\,
\frob^{\boundm *}(\shS)$ is semistable, then $\shS$ is strongly
semistable.
\end{lemma}
\proof By enlarging the ground field we may assume that there exists
an $\fifi$-rational point on $\cu$ and hence an invertible sheaf
$\shL$ of degree $1$. Let $\shS$ be given of rank $r$. There exists
$a \geq 0$ and $e >0$, $a+e \leq r$, such that $p^{a+e} \deg(\shS)
=p^a \deg(\shS) \modu r$. Let $\shT= \frob^{a*}(\shS) \otimes
\shL^\ell$, where $\ell$ is chosen such that $\delta := \deg(\shT)=
p^a \deg(\shS) + \ell r$ is in between $0 $ and $r-1$. There exists
$k$ such that $p^{e} \delta - \delta = kr$. Hence the assignment
(and its iterations)
$$\shT \longmapsto \frob^{e*}(\shT) \otimes\shL^{-k} =: \varphi(\shT)$$ preserves the degree, since $\deg
(\frob^{e*}(\shT) \otimes \shL^{-k})= p^{e} \deg (\shT) -kr$. The
number of semistable bundles defined over $\fifi$ of fixed rank $r$
and degree $\delta$ is finite, say bounded by $\boundn$. Set
$\boundm = r (\boundn +1)$.

If now $\frob^{\boundm*} (\shS)$ is semistable, then all
$\frob^{i*}(\shS)$ are semistable for $i \leq \boundm $, but then of
course the $\varphi^{j} (\shT)$ are semistable for all $j \leq
\boundn$. Since they have the same degree $\delta$, two of them must
be isomorphic. But then we must have a periodicity among the
$\varphi^{j}(\shT)$ and so they are in fact semistable for all $j
\in \NN$. Hence $\shS$ is strongly semistable (see also
\cite{langestuhler} for this argument).

Note that we have to deal only with the finitely many degrees
between $0$ and $r-1$, so that there exists a bound which is
independent of the degree of the bundles. \qed

\medskip
For every locally free sheaf $\shS$ on $\cu$ there exists the
so-called Harder-Nara\-sim\-han filtration $0 \subset \shS_1
\subsetdots \shS_t =\shS$, where the $\shS_j$ are locally free
subsheaves. This filtration is unique and has the property that the
quotients $\shS_{j}/ \shS_{j-1}$ are semistable and
$\mu(\shS_{j}/\shS_{j-1})>\mu(\shS_{j+1}/ \shS_{j})$ for all $j= 1
\komdots t-1$. In positive characteristic, a Harder-Narasimhan
filtration is called strong if all quotients $\shS_j/\shS_{j-1}$ are
strongly semistable. For every locally free sheaf $\shS$ there
exists a pull-back $\frob^{e*} (\shS)$ such that its
Harder-Narasimhan filtration is strong \cite[Theorem
2.7]{langersemistable}. An observation of V. Trivedi combined with
Lemma \ref{strongtest} allows us to give a bound for the Frobenius
power such that the Harder-Narasimhan filtration is strong.

\begin{corollary}
\label{hnbound} Let $\cu$ denote a smooth projective curve of genus
$g$ defined over a finite field $\fifi$, and let $r \in \NN$.
Suppose that $\Char (\fifi) = p
>4(g-1) r^3$.
Then there exists a number $\boundhn$ such that for every locally
free sheaf $\shS$ of rank $r$ the Harder-Narasimhan filtration of
$\, \frob^{\boundhn*}( \shS)$ is strong.
\end{corollary}
\proof We use induction on $r$; for $r=1$ there is nothing to show.
Set $\boundhn = (r-1) \boundm$, where $\boundm$ is the number for
which the conclusion of Lemma \ref{strongtest} holds for all locally
free sheaves of rank $\leq r$. Consider $\frob^{\boundm *}(\shS)$.
If this is semistable, then $\shS$ is strongly semistable by Lemma
\ref{strongtest} and we have the desired result. So suppose that
$$0 \subset \shS_1
\subsetdots \shS_{t-1} \subset \shS_t= \frob^{\boundm *} (\shS)
$$ is the Harder-Narasimhan filtration, with $\rk (\shS_j/ \shS_{j-1})
<r$. By the Theorem of Trivedi \cite{trivedihilbertkunzreduction}
(here the condition about the prime characteristic is needed) the
Harder-Narasimhan filtration of a higher Frobenius pull-back is a
refinement of the pull-back of this filtration. Hence we can apply
the induction hypothesis to the quotient sheaves $\shS_j/
\shS_{j-1}$. \qed

\begin{remark}
For $r=2$ we do not need the condition about the prime
characteristic, since then $\frob^{\boundm}(\shS)$ is either
semistable ($\boundm$ being the invariant from Lemma
\ref{strongtest}), and hence strongly semistable, or it has an
invertible subsheaf which contradicts semistability. In any case the
Harder-Narasimhan filtration of $\frob^{\boundm}(\shS)$ is strong.
\end{remark}

\section{Frobenius test exponents for curves of higher genus}

Let $\cu$ denote a smooth projective curve over a finite field and
let $\shS$ denote a locally free sheaf on $\cu$. We say that a
cohomology class $c \in H^1(\cu,\shS)$ is annihilated by some
Frobenius if $\frob^{e*} (c) \in H^1(\cu, \frob^{e*}(\shS))$ is zero
for some $e \in \NN$. We want to show that this can be checked
within a certain number of steps, which is bounded by a number only
depending on $\cu$ and on the rank of $\shS$. We restrict first to
strongly semistable sheaves.

\begin{lemma}
\label{strongsemistabletestexponent} Let $\cu$ denote a smooth
projective curve defined over a finite field $\fifi$, and let $r \in
\NN$. Then there exists a number $\boundc$ such that for every
strongly semistable locally free sheaf $\sheS$ on $\cu$ of rank $r$
the Frobenius annihilation of a cohomology class $c \in H^1(\cu,
\sheS)$ can be checked within $\boundc$ steps.
\end{lemma}
\proof We examine three cases depending on whether $\deg (\sheS)$ is
positive, zero or negative. In each case we will find a bound
$\boundc$.

Suppose first that $\deg (\sheS) >0$. Let $\boundc$ be such that
$p^{\boundc} > r \deg (\omega_{\cu})$, where $\omega_{\cu}$ denotes
the dualizing sheaf on $\cu$; thus $\deg(\omega_{\cu})= 2g-2$, where
$g$ is the genus of the curve. Then $\frob^{\boundc*}(c) \in
H^1(\cu, \frob^{\boundc *}(\sheS))$ and the degree of this bundle is
$$ \deg (\frob^{ \boundc *} (\sheS)) = p^{\boundc} \deg(\sheS) \geq p^\boundc \, . $$
Therefore $\mu ( \frob^{ \boundc *} (\sheS)) \geq p^\boundc/ r >
\deg (\omega_{\cu})$. Hence there exists no non-trivial homomorphism
$\frob^{ \boundc *} (\sheS) \ra \omega_{\cu}$ (due to strong
semistability) and therefore $H^1(\cu, \frob^{ \boundc *}
(\sheS))=0$ by Serre duality. So the $\boundc$-th Frobenius
annihilates these classes.

Suppose now that $\deg(\sheS)=0$. The curve and $\sheS$ are defined
over the finite field $\fifi$. The pull-back of $\sheS$ is again
semistable of degree $0$ and defined over $\fifi$. The number of
semistable bundles of rank $r$ and degree $0$ defined over a finite
field is however finite. Call this number $n$, as in Lemma
\ref{strongtest}. Hence there must be a repetition, say $\frob^{t *}
(\sheS) \cong \frob^{t' *}(\sheS)$, where $0 \leq t <t' \leq n$. Set
$t'-t =: v$ and define $\shF:= \frob^{t *}(\sheS)$.

Let $\theta : \frob^{v*}(\shF) \ra \shF$ be a fixed isomorphism.
This induces isomorphisms $\frob^{(k-1)v*}(\theta):
\frob^{kv*}(\shF) \ra \frob ^{(k-1)v*}(\shF)$ and the composition of
these yields isomorphisms $\theta_k : \frob^{kv*}(\shF) \ra \shF$.
Set $\varphi = \theta \circ \frob^{v*} :\shF \ra \shF$. The diagram
$$
\begin{CD}
\frob^{v*}(\shF) @> \frob^{v*} >> \frob^{v*} (\frob^{v*}(\shF)) \cr
@V \theta VV         @VV \frob^{v*}( \theta) V \cr
\shF@>\frob^{v*}>>\frob^{v*}(\shF)
\end{CD}
$$
commutes by the functoriality of $\frob^*$. Therefore $\theta_k
\circ \frob^{kv*} = \varphi^k$. Since the $\theta_k$ are
isomorphisms, they induce isomorphisms on the cohomology groups, and
so a cohomology class $c \in H^1(\cu, \shF)$ is annihilated by
$\frob^{kv*}$ if and only if it is annihilated by $ \varphi^k$.

So we are considering the group homomorphism $\varphi: H^1(\cu,
\shF) \ra H^1(\cu, \shF)$ and its iterations. Note that $\varphi$ is
$p^v$-linear, meaning that $\varphi(\lambda c)=
\lambda^{p^v}\varphi(c)$. In this situation we have an
$\fifi$-vector space decomposition $H^1(\cu, \shF) =V_s \oplus V_n$
such that the action of $\varphi$ on $V_s$ is bijective and on $V_n$
is nilpotent \cite[14, last corollary]{mumfordabelian}. Now a
cohomology class $c \in H^1(\cu,\sheS)$ is annihilated by some
Frobenius power if and only if the image of $c$ under ${\rm proj}_s
\circ \frob^{t*}: H^1(\cu, \sheS) \ra H^1(\cu, \shF) \ra V_s$ is
zero. The dimension of $H^1(\cu, \shF)$ is bounded by an invariant
dependent only on $r$ and $\cu$, which follows from \cite[Lemma
20]{seshadrifibre}. Every application of $\varphi$ annihilates at
least one dimension of $V_n$, so there is a certain power of
$\varphi$ which annihilates $V_n$.

Now suppose that $\deg (\sheS) < 0$. Even in this case it might
happen that some Frobenius power of a non-zero element $c \in
H^1(\cu, \sheS)$ is zero, but this is rather an exception. First
note that the dual sheaf $\sheS^\dual$ has positive degree and is
therefore ample, since it is strongly semistable. Consider the
extension $0 \ra \sheS \ra \sheS' \ra \O_{\cu} \ra 0$ defined by $c
\in H^1(\cu,\sheS)$ and consider the dual sequence $0 \ra \O_{\cu}
\ra (\sheS')^\dual \ra \sheS^\dual \ra 0$ and its Frobenius
pull-backs $$0 \lra \O_{\cu} \lra \frob^{t*}(\sheS')^\dual \lra
\frob^{t*}(\sheS^\dual) \lra 0 \, .$$ If, for some $t$, this
extension is still non-trivial, it follows from \cite[Proposition
2.2]{giesekerample} that every quotient bundle of
$\frob^{t*}(\sheS')^\dual$ has positive degree. Now let $\boundr$ be
such that
$$ \deg ( \frob^{ \boundr *} (\sheS')^\dual) = \deg (\frob^{ \boundr *} (\sheS)^\dual)
=p^\boundr  \deg(\sheS^\dual) >\rk (\sheS) g$$ (note that this is
satisfied if $p^\boundr>\rk(\sheS)g$). Then either
$\frob^{\boundr*}(c)=0$ and $c$ is annihilated by this Frobenius
power, or $\frob^{\boundr*}(c) \neq 0$, the extension remains
non-trivial and we have, by \cite[Lemma 2.2]{giesekerample}, that $
\frob^{ \boundr *} (\sheS')^\dual$ is ample. But then the torsor
defined by $c$, that is $\PP((\shT') ^\dual) - \PP(\shT^\dual)$, is
an affine variety. Therefore $c$ cannot be annihilated by any
Frobenius power, because otherwise the pull-back of this torsor
would admit a section, contrary to affineness. \qed

\begin{theorem}
\label{testexponentcurve} Let $\cu$ denote a smooth projective curve
of genus $g$ defined over a finite field $\fifi$. Fix $r \in \NN$
and suppose that $\Char (\fifi)= p>4(g-1) r^3$. Then there exists a
number $\bounda$ such that for every locally free sheaf $\shS$ of
rank $r$ on $\cu$ the following holds: if a cohomology class $c \in
H^1(\cu, \shS)$ is annihilated by some Frobenius power, then already
$\frob^{\bounda *}(c) = 0$ holds in $H^1(\cu, \frob^{\bounda *}
(\shS))$.
\end{theorem}
\proof By Corollary \ref{hnbound} there exists a Frobenius power
$\frob^{\boundhn}$ such that the Harder-Narasimhan filtration of
$\frob^{\boundhn*}(\shS)$ is strong, say
$0\subset\shS_1\subsetdots\shS_t = \frob^{\boundhn*}(\shS)$. Let
$\shT = \shS_j/\shS_{j-1}$ be one of the strongly semistable
quotient sheaves in this filtration, so that $\rk (\shT) \leq r$. By
Lemma \ref{strongsemistabletestexponent} we know that for all these
bundles there exists a checking bound for Frobenius annihilation.
Let $\boundf$ be such a common bound and set $\bounda = \boundhn+r
\boundf$.

Now let $c \in H^1(\cu, \shS)$ denote a cohomology class and suppose
that it is annihilated by some Frobenius power. We want to show that
already the $\bounda$-th Frobenius power annihilates this class. Let
$c_t= \frob^{\boundhn*}(c) \in H^1(\cu, \frob^{\boundhn *}(\shS))$
be the pull-back and consider what happens to it under the short
exact sequence
$$ 0 \lra \shS_{t-1} \lra \frob^{\boundhn *}(\shS) \lra
\frob^{\boundhn *}(\shS)/ \shS_{t-1} \lra 0 \, .$$ Since some
Frobenius power of $c_t$ is zero, this is also true for its image
$c_t' \in H^1(\cu,\frob^{\boundhn *}(\shS)/ \shS_{t-1} )$. But then
$\frob^{\boundf*}(c_t') =0 $ in $H^1(\cu, \frob^{ \boundf *}(
\frob^{\boundhn *}(\shS)/ \shS_{t-1}))$ by Lemma
\ref{strongsemistabletestexponent}. Therefore $\frob^{\boundf*}
(c_t)$ is the image of a cohomology class $c_{t-1} \in H^1(\cu,
\frob^{\boundf *} (\shS_{t-1}) )$.

If $\shS_t/\shS_{t-1}$ has negative degree, then
$H^0(\cu,\shS_t/\shS_{t-1})= 0$ by semistability and $H^1(\cu,
\shS_{t-1}) \ra H^1(\cu, \shS_t)$ is injective, and this holds for
all Frobenius powers. Hence also $c_{t-1}$ is annihilated by some
Frobenius power. If however $\deg(\shS_t/\shS_{t-1}) \geq 0$, then
$\mu_{\rm min} (\shS_{t-1})= \mu (\shS_{t-1}/ \shS_{t-2})
> 0$ and every cohomology class in it is annihilated by some
Frobenius power. Application of the induction hypothesis to
$\shS_{t-1}$ gives the result. \qed

\begin{remark}
The bound obtained in Theorem \ref{testexponentcurve} is hardly
suitable for computations. The main problem here is the number $n$
from Lemma \ref{strongtest}, which bounds the number of semistable
sheaves of given rank and degree on a curve over a finite field
(this number enters also in the second case in the proof of Lemma
\ref{strongsemistabletestexponent}). The dimension of the moduli
space of semistable sheaves of rank $r$ and degree $d$ is
$r^2(g-1)+1$ (see \cite[after Th\'{e}or\`{e}me 18]{seshadrifibre} or
\cite[Corollary 4.5.5]{huybrechtslehn}; note that Huybrechts and
Lehn give the dimension for a fixed determinant, and this explains
the difference of $g$).
\end{remark}

\section{Frobenius test exponents for ideals in two-dimensional rings}

We come now back to Frobenius test exponents for ideals in
two-dimensional rings.

\begin{theorem}
\label{testhighergenus} Let $\fifi$ denote a finite field and let
$R$ denote a geometrically normal two-dimensional standard-graded
domain over $\fifi$. Fix $n \in \NN$ and suppose that $\Char
(\fifi)=p \geq 4g(n-1)^3$, where $g$ denotes the genus of the smooth
projective curve $\cu= \Proj R$. Then there exists a Frobenius test
exponent for the class of homogeneous ideals generated by at most
$n$ homogeneous elements.
\end{theorem}
\proof Suppose first that $I$ is an $R_+$-primary ideal. Let $I=(f_1
\komdots f_n)$ and suppose as in the proof of Theorem
\ref{testelliptic} that $f \in I^F$ is homogeneous of degree $m$.
Let $\delta(f) \in H^1(\cu,\Syz(f_1 \komdots f_n)(m))$. Then
$f^{p^e} \in I^{[p^e]}$ if and only if $\delta(f)$ is annihilated by
the $e$-th Frobenius power. By Theorem \ref{testexponentcurve}, we
know that $\delta(f)$ is annihilated already by the $\bounda$-th
Frobenius power, where $\bounda$ is the test bound for locally free
sheaves of rank $r=(n-1)$.

Suppose now that $I$ is a homogeneous ideal, but not necessarily
$R_+$-primary, say $I=(h_1 \komdots h_k)$. Suppose that $f$ is a
homogeneous element of degree $m$. Let $x,y$ denote homogeneous
parameters in $R$ of degree $> m$. If now $f \in I^F$, then also $f
\in (I +(x,y))^F$ (which is the Frobenius closure of an
$R_+$-primary ideal generated by $k+2$ elements), and hence
$f^{p^\boundb} \in (h_1^{p^\boundb} \komdots
h_k^{p^\boundb},x^{p^\boundb},y^{p^\boundb})$ (where $\boundb$ is
now the bound for sheaves of rank $k+1$). If we write this as a
homogeneous equation $f^{p^\boundb} = g_1 h_1^{p^{\boundb}}
\plusdots g_k h_k^{p^{\boundb}} + g_{k+1}x^{p^\boundb} + g_{k+2}
y^{p^\boundb}$ we see that already $f^{p^\boundb} \in
I^{[p^\boundb]}$. \qed

\begin{remark}
If we delete the condition on the prime number in Theorem
\ref{testhighergenus}, we still get a positive answer to the weak
question of Katzman and Sharp for primary ideals in a
two-dimensional geometrically normal standard-graded domain over a
finite field. If $I=(f_1 \komdots f_n)$ is the primary ideal, then
its syzygy bundle has a strong Harder-Narasimhan filtration on a
certain Frobenius pull-back: let us say that the Harder-Narasimhan
filtration of the $\boundhn$-th pull-back is strong. Now the syzygy
bundle of $I^{[q]}=(f_1^q \komdots f_n^q)$, where $q=p^{e}$, that is
$\Syz(f_1^q \komdots f_n^q)$, is just the $e$-th pull-back of
$\Syz(f_1 \komdots f_n)$, and therefore the $\boundhn$-th pull-back
of $\Syz(f_1^q \komdots f_n^q)$ has also a strong Harder-Narasimhan
filtration. This replaces Corollary \ref{hnbound}.
\end{remark}

\section{Tight closure test ideal exponents}
\label{tightexponent}

We recall briefly the notions of tight closure, test ideals and test
exponent for tight closure, referring to \cite{hunekeparameter} and
\cite{hochsterhuneketestexponent} for details. Let $R$ denote a
noetherian commutative ring of positive prime characteristic $p$,
and let $R^\minpricomp$ denote the complement of the union of the
minimal prime ideals. Then the tight closure of an ideal $I$ is
defined as the ideal
$$I^*= \{f \in R:\, \exists \, z \in R^\minpricomp \mbox{
such that } zf^q \in I^{[q]} \mbox{ for all } q \gg 0 \} \, .$$

An element $z \in R^\minpricomp$ is called a test element if for all
ideals $I$ and all $f \in I^*$ we have $zf^q \in I^{[q]}$ for all
powers $q=p^{e}$. The test ideal, denoted $\test$, is the ideal
generated by all test elements. It is a non-trivial fact that test
elements exist \cite[Theorem 3.2]{hunekeparameter} under certain
conditions. The situation for Gorenstein local or graded rings is
even better: in this case for an arbitrary system of parameters we
have $(x_1 \komdots x_d): \test = (x_1 \komdots x_d)^*$ and
$\test=(x_1 \komdots x_d):(x_1 \komdots x_d)^*$; see
\cite[Corollaries 4.2 and 4.3]{hunekeparameter}. More specifically,
for a noetherian nonnegatively graded Gorenstein ring and $p \gg 0$
we have that $\test = R_{\geq a+1}$, where $a$ is the a-invariant of
$R$. This is the maximal degree $\delta$ such that
$(H^d_{R_+}(R))_\delta \neq 0$; equivalently, $\O_{\cu}(a)$ is the
dualizing sheaf on $\cu = \Proj R$ \cite[Proposition
3.6.11]{brunsherzog}.

Fix a test element $z \in R$, and let $I \subseteq R$ be an ideal. A
test exponent for $z$ and $I$ is a prime power $p^{b}$ (or rather
its exponent) such that $zf^{p^b} \in I^{[p^b]}$ ensures that $f \in
I^*$ (see \cite[Definition 2.2]{hochsterhuneketestexponent}, where
the definition is given for $R$-modules $N \subseteq M$).

Test exponents are important for two reasons. On the one hand, the
existence of test exponents is equivalent to the localization of
tight closure; see \cite[Theorem 2.4]{hochsterhuneketestexponent}
for the precise statement. On the other hand, a test exponent gives
at once a finite algorithm for the computation of tight closure.
Note that \--- contrary to the case of Frobenius test exponents \---
already the existence of a test exponent for a single ideal is a
problem, not to mention the existence of a uniform bound for test
exponents for a reasonable class of ideals \cite[Discussion
5.3]{hochsterhuneketestexponent}.

Here we will however focus on the following variant of test
exponent, which is easier to handle.

\begin{definition}
Let $R$ denote a noetherian ring of positive characteristic $p$ and
let $\test$ denote the test ideal. A test ideal exponent for tight
closure (for a certain class of ideals in $R$) is a number
$\boundtight$ such that the following holds: if $z f^{p^e} \in
I^{[p^e]}$ holds for every $e \leq \boundtight$ and for every $z\in
\test $, then $f \in I^*$.
\end{definition}

The existence of a test ideal exponent has the same computational
impact on tight closure as the existence of a test exponent for a
certain test element. Moreover, its existence implies also that
tight closure commutes with localization, at least if all the test
elements are locally stable (look at the proof of \cite[Proposition
2.3]{hochsterhuneketestexponent}).

We restrict now to the situation of a two-dimensional
standard-graded normal Gorenstein domain $R$. We will show that in
this case a uniform tight closure test ideal exponent exists for the
ideals generated by $n$ homogeneous elements. For a homogeneous
element $f$ of degree $m$ and a homogeneous element $z$ the
condition that $zf^q \in (f_1 \komdots f_n)^{[q]}$ holds translates
into $z\frob^{e*} (\delta(f))=0$ in
$H^1(\cu,\O_{\cu}(\deg(z))\otimes\frob^{e*}(\shS))$, where $\shS =
\Syz(f_1 \komdots f_n)(m)$, $\delta (f) \in H^1(\cu, \shS)$ and
$q=p^{e}$. We will need the following lemma.

\begin{lemma}
\label{negative} Let $\cu$ denote a smooth projective curve of genus
$g$ over a field $\fifi$. Fix a very ample invertible sheaf
$\O_{\cu}(1)$ and set $\deg(\cu) =\deg (\O_{\cu}(1))$. Let $\bonum$
be a number such that $\bonum \deg(Y) > 2g-2= \deg(\omega_{\cu})$,
where $\omega_{\cu}$ denotes the dualizing sheaf on $\cu$. Then for
every semistable sheaf $\shT$ on $\cu$ of rank $r$ and of degree $
\deg (\shT) < -r(1+(\bonum+1)\deg(Y))$ and for every cohomology
class $0 \neq c \in H^1(\cu, \shT)$ there exists a sheaf
homomorphism $\varphi: \shT \ra \O_{\cu}(- \bonum ) \otimes
\omega_{\cu}$ such that $\varphi(c) \neq 0$ in $ H^1(\cu,\O_{\cu}(-
\bonum ) \otimes \omega_{\cu} )$.
\end{lemma}
\proof Suppose that $\shT$ satisfies the stated degree condition.
Then
\begin{eqnarray*} \deg(\shT^\dual
\otimes \omega_{\cu} \otimes \O_{\cu}(-\bonum -1)) & =& - \deg(\shT)
+ r(2g-2)-r(\bonum+1)\deg(Y) \cr
&>& r (2g-1) \, .
\end{eqnarray*} By \cite[Lemme 20]{seshadrifibre} the sheaf
$\shT^\dual \otimes \omega_{\cu} \otimes \O_{\cu}(- \bonum-1)$ is
therefore generated by its global sections and $H^1(\cu,\shT^\dual
\otimes \omega_{\cu} \otimes \O_{\cu}(-\bonum -1) )=0$. This last
property means that $\shT^\dual \otimes \omega_{\cu} \otimes
\O_{\cu}(- \bonum )$ is $0$-regular in the sense of
Castelnuovo-Mumford \cite[Definition 1.7.1]{huybrechtslehn}.
Therefore by \cite[Lemma 1.7.2]{huybrechtslehn} we have a surjective
map
$$H^0(\cu, \shT^\dual \otimes \omega_{\cu} \otimes \O_{\cu}(-\bonum))\otimes H^0(\cu, \O_{\cu}(n))
\lra H^0(\cu,\shT^\dual \otimes \omega_{\cu} \otimes\O_\cu(n-\bonum)
) \,
$$ for every $n \geq 0$. Taking $n= \bonum$ we get a surjective
mapping $\varphi:\bigoplus \O_{\cu}(\bonum) \ra  \shT^\dual \otimes
\omega_{\cu} $ which is also globally surjective. Consider the exact
sequence
$$ 0 \lra \shK
\lra \bigoplus \O_{\cu}(\bonum) \stackrel{\varphi}{\lra} \shT^\dual
\otimes \omega_{\cu} \lra 0 \, .$$ This induces the exact sequence
$$ \bigoplus H^0(\cu, \O_{\cu}(\bonum)) \stackrel{\varphi}{\lra} H^0(\cu,
\shT^\dual \otimes \omega_{\cu}) \lra H^1(\cu, \shK) \lra \bigoplus
H^1(\cu, \O_{\cu}(\bonum)) \, .$$ Since $H^1(\cu, \O_{\cu}(\bonum))=
H^0(C , \O_{\cu}(- \bonum) \otimes \omega_{\cu} )^\dual =0$ in view
of the assumption on $u$, and since $\varphi $ is surjective on the
global sections by construction, we get $H^1(\cu, \shK)=0$. Now
tensor the short exact sequence with $\omega_{\cu}^{-1}$ and
consider the dual sequence
$$0 \lra \shT \lra \bigoplus \O(-\bonum) \otimes \omega_{\cu} \lra \shK^\dual
\otimes \omega_{\cu} \lra 0 \, .$$ Since $H^0(\cu, \shK^\dual
\otimes \omega_{\cu})= H^1(\cu, \shK)^\dual =0$, we see that the
induced map
$$H^1(\cu,\shT) \lra \bigoplus H^1(\cu, \O_{\cu}(-\bonum) \otimes \omega_{\cu}))$$ is
injective. Hence for a class $ 0 \neq c \in H^1(\cu, \shT)$ at least
one component in $H^1(\cu, \O_{\cu}(-\bonum) \otimes \omega_{\cu}))$
is non-zero. \qed


\begin{lemma}
\label{negativefrobenius} Let $R$ denote a standard-graded
two-dimensional geometrically normal Gorenstein domain over a field
$\fifi$ of positive characteristic $p$. Let $\cu = \Proj \, R$
denote the smooth projective curve of genus $g$ determined by $R$.
Set $\deg(Y)= \deg(\O_{\cu}(1))$ and let $\omega_{\cu}=\O_{\cu}(a)$
be the dualizing sheaf. Let $\test $ denote the test ideal of $R$.
Fix $r \in \NN$. Let $\boundo$ be such that $2^\boundo > r(
1+(a+g+2) \deg(Y))$. Then $\boundo$ has the property that for every
strongly semistable sheaf $\shT$ of rank $r$ and of negative degree
the following holds: if for $c \in H^1(\cu, \shT)$ we have that $z
\frob^{\boundo*} (c) = 0$ for all homogeneous $z \in \test $, then
already $ \frob^{\boundo *} (c)=0$.
\end{lemma}
\proof Since $\shT$ has negative degree, the degree of the
$\boundo$-th Frobenius pull-back $\frob^{ \boundo *}(\shT)$ is
$$\deg(\frob^{\boundo *} (\shT)) =p^{\boundo} \deg (\shT) \leq
-2^\boundo < -r( 1+(a+g+1) \deg(Y)) \, .$$ We want to apply Lemma
\ref{negative} with $\bonum = a+g+1$; note that $(a+g+1) \deg(\cu) >
\deg \O_\cu(a) = 2g-2$. Suppose that $\frob^{\boundo *}(c) \neq 0$.
Then there exists by Lemma \ref{negative} a sheaf homomorphism
$\varphi: \frob^{\boundo*}(\shT) \ra \O_{\cu}(a- \bonum)=
\O_{\cu}(-g-1)$ such that $0 \neq c'=\varphi(\frob^{\boundo*}(c))
\in H^1(\cu, \O_{\cu}(-g-1))$. Since $z \frob^{\boundo *}(c)=0$ for
all $z \in \test$, this holds also for $c'$. This means that we have
a cohomology class $c' \in H^1(\cu, \O(-g-1 )) \cong
(H^2_{R_+}(R))_{-g-1}$ of negative degree $-g-1$ which is
annihilated by the test ideal. By \cite[Corollary
4.2(1)]{hunekeparameter} we have that $0^* = \{ c \in
H^2_{R_+}(R):\test c = 0 \}$ and hence $c' \in 0^*$. This means that
the torsor defined by $c'$ is not affine, but since
$\deg(\O_{\cu}(-g-1))= -(g+1) \deg(\cu) < -g$ the theorem of
Gieseker (\cite[Lemma 2.2]{giesekerample}, see the third case of the
proof of Lemma \ref{strongsemistabletestexponent} above) implies
that $c'=0$, a contradiction. \qed

\begin{theorem}
\label{tightexponenttheorem} Let $R$ denote a standard-graded
two-dimensional geometrically normal Gorenstein domain over a finite
field $\fifi$ of positive characteristic $p$. Fix $n \in \NN$ and
suppose that $p \geq 4 g (n-1)^3 $, where $g$ denotes the genus of
$\cu = \Proj R$. Then there exists a test ideal exponent for the
class of homogeneous ideals generated by at most $n$ elements.
\end{theorem}
\proof We have to show that there exists a $\boundtight$ such that
for every homogeneous ideal $I=(f_1 \komdots f_n)$ the following
holds: if for an element $f \in R$ we have $zf^{q} \in I^{[q]}$ for
all $z \in \test$ and all $q=p^{e}$, $e \leq \boundtight$, then
already $f \in I^*$.

Again we may reduce to the primary case as in the proof of Theorem
\ref{testhighergenus}. Set $\boundtight=\boundhn+(n-1)\boundo $,
where $\boundhn$ is the bound from Corollary \ref{hnbound} for
locally free sheaves of rank $n-1$ on $\cu$ and where $\boundo$ is a
common bound from Lemma \ref{negativefrobenius} for the strongly
semistable sheaves of ranks $r \leq n-1$. We can also reduce to the
homogeneous case, so let $f$ denote a homogeneous element of degree
$m$. Suppose that $zf^q \in I^{[q]}$ for all $z \in \test$ and all
$q=p^e$, $e \leq \boundtight$; we have to show that this is then
true for every power of $p$ (and at least for one element $z \in
R^\minpricomp$). Let $\shS= \Syz(f_1 \komdots f_n)(m)$ denote the
syzygy bundle of rank $n-1$ on the smooth projective curve $\cu =
\Proj R$ and let $c=\delta(f) \in H^1(\cu, \Syz(f_1 \komdots
f_n)(m))$.

Using Corollary \ref{hnbound} we get the strong Harder-Narasimhan
filtration $\shS_0 \subsetdots \shS_t =\frob^{\boundhn*}(\shS)$;
consider $c_t = \frob^{\boundhn *}(c) \in
H^1(\cu,\frob^{\boundhn*}(\shS))$. Then we consider the image $c_t'$
of $c_t$ in $H^1( \cu, \shS_t /\shS_{t-1})$. If $\shT = \shS_t
/\shS_{t-1}$ has negative degree, then by Lemma
\ref{negativefrobenius} we have that $\frob^{\boundo*}(c_t')=0$,
since it is annihilated by the test ideal. In this case
$\frob^{\boundo *}(c_t)$ is the image of a class $c_{t-1} \in
H^1(\cu, \frob^{\boundo *}(\shS_{t-1}))$. Since $H^0(\cu,\frob^{e
*}( \shT))=0$, we have $H^1(\cu,\frob^{e*}(\shS_{t-1}))\subseteq
H^1(\cu,\frob^{e*}(\shS_t))$ for all $e$ and it follows from our
assumption that $z \frob^{e*}(c_{t-1})=0$ for $e\leq (n-2) \boundo $
and all $z \in \test$. If however $\shT$ has non-negative degree,
then the minimal slope of $\shS_t=\frob^{\boundhn*}(\shS)$ is
non-negative and then $z \frob^{e* }(c) = 0$ for $z \in R_{\geq
a+1}$  and all $e$ anyway (even in small characteristics, where
$\test$ may be different from $R_{\geq a+1}$). In this way we use
induction along the strong Harder-Narasimhan filtration. \qed

\begin{remark}
The problem with test exponents for a given fixed test element $z$
is that even in the case of a parameter ideal (which corresponds to
a cohomology class in an invertible sheaf) it is not clear how to
bound the test exponent \cite[Discussion
5.1]{hochsterhuneketestexponent}.
\end{remark}

\begin{corollary}
\label{testcubic} Let $\fifi$ denote a field of positive
characteristic $p$ and let $G \in \fifi[x,y,z]$ denote a cubic
polynomial such that $\cu= \Proj R$, where $R=\fifi[x,y,z]/(G)$, is
an elliptic curve. Let $e$ be such that $p^{e} >7 (n-1)$. Then $e$
is a test ideal exponent for the set of all homogeneous
$R_+$-primary ideals generated by $n$ homogeneous elements.
\end{corollary}
\proof Since $g(\cu)=1$, $\deg (\cu)=3$ and $\omega_\cu = \O_\cu $,
Lemma \ref{negative} applied with $u=1$ shows that for every
semistable sheaf $\shT$ on $\cu$ of rank $r$ and of degree $\deg
(\shT) < -7r$ and a cohomology class $c \in H^1(\cu,\shT)$ there
exists a sheaf homomorphism $\varphi: \shT \ra \O_\cu(-1)$ such that
$\varphi(c) \neq 0$ in $H^1(\cu, \O_\cu(-1))$.

Let $q=p^{e} >7(n-1)$. Let $I=(f_1 \komdots f_n)$ denote an
$R_+$-primary ideal generated by $n$ homogeneous elements and let $f
\in R$ denote an element of degree $ m$. The test ideal of $R$ is
the maximal ideal $\test=(x,y,z)$. So suppose that $wf^q \in (f_1^q
\komdots f_n^q)=I^{[q]}$ for $w=x,y,z$. We have to show that $f \in
I^*$. Let $\Syz(f_1 \komdots f_n)(m) \cong \shS_1 \oplusdots \shS_k$
denote the decomposition of the syzygy bundle on $\cu$ into
indecomposable sheaves and let $$\delta(f)=c=(c_1 \komdots c_k) \in
H^1(C, \Syz(f_1 \komdots f_n)(m))= \bigoplus_{j=1}^k
H^1(\cu,\shS_j)$$ denote the components of the cohomology class.
Assume that $f \not\in I^*$. By \cite[Corollary
4.1]{brennertightelliptic} there exists $j$, $1 \leq j \leq k$, such
that $\shS_j$ has negative degree and $c_j \neq 0$. Consider
$\frob^{e*}(c) \in H^1(\cu ,\frob^{e*}(\shS_j)$). Recall that
indecomposable sheaves on an elliptic curve are strongly semistable.
Furthermore note that $\frob^{e*}(c) \neq 0$, otherwise $f$ would
belong to the tight closure. Since $\deg( \frob^{e*}(\shS_j)) \leq
p^{e} < -7(n-1) \leq -7 \rk(\shS_j)$, there exists a morphism
$\varphi: \frob^{e*}(\shS_j) \ra \O_\cu(-1) $ such that $0\neq
c'=\varphi(\frob^{e*}(c_j))\in H^1(\cu, \O_C(-1))$. By assumption
the class $\frob^{e*}(c_j)$ is annihilated by $(x,y,z)$, and this
property passes over to $c'$. But for $0 \neq c' \in H^1(\cu,
\O_\cu(-1))$ there exists, by Serre duality, a homomorphism $\theta
: \O_C(-1) \ra \O_\cu = \omega_\cu$ such that $\theta(c') \neq 0$.
However, $\theta$ is given by a linear form, and so we have a
contradiction. \qed

\bibliographystyle{plain}

\bibliography{bibliothek}

\end{document}